\documentclass[a4paper,11pt]{amsart}
\usepackage[latin1]{inputenc}				
\usepackage[T1]{fontenc}					
\usepackage{lmodern}						
\usepackage[english]{babel}					
\usepackage{amsmath, latexsym, amsfonts, amssymb, amsthm, amscd}			
\usepackage{amsmath}
\usepackage{pgfplots}
\usepackage{tikz}
\usetikzlibrary{calc}
\usetikzlibrary{arrows,backgrounds}
\tikzstyle{NN}=[circle,fill=black,draw=black,
inner sep=0pt,minimum size=2mm]
\tikzstyle{NR}=[circle,fill=white,draw=black,
inner sep=0pt,minimum size=2mm]
\usepackage{graphicx}						
 \newenvironment{Proof}{
  \removelastskip\par\medskip  
  \noindent{} {\bf Proof.}\rm}{\penalty-20\null\hfill$\square$\par\medbreak}

\newcommand{\ga}{\alpha}
\newcommand{\gb}{\beta}
\newcommand{\gga}{\gamma}            % \gg already exists...

\newcommand{\gd}{\delta}

\newcommand{\gl}{\lambda}

\newcommand{\gs}{\sigma}

\newcommand{\cE}{{\ensuremath{\mathcal E}} }

\newcommand{\cO}{{\ensuremath{\mathcal O}} }

\newcommand{\cU}{{\ensuremath{\mathcal U}} }

\newcommand{\C}{{\ensuremath{\mathbb C}} }

\newcommand{\E}{{\ensuremath{\mathbb E}} }

\newcommand{\bbP}{{\ensuremath{\mathbb P}} }

\newcommand{\R}{{\ensuremath{\mathbb R}} }

\newcommand{\ovl}{\overline}
\newcommand{\rrb}{\right\rbrace}
\newcommand{\llb}{\left\lbrace}

\newcommand{\Cov}{\mathrm{Cov}}
\newcommand{\dis}{\mathrm{distribution}}

\newtheorem{theorem}{Theorem}[section]
\newtheorem{lemma}[theorem]{Lemma}
\newtheorem{prop}[theorem]{Proposition}

\newtheorem{rem}[theorem]{Remark}
\newtheorem{definition}[theorem]{Definition}

\newtheorem{hyp}[theorem]{Hypothesis}
\title[Eigenvectors of Sample Covariance Matrices]{Eigenvectors of Sample Covariance Matrices: Universality of global fluctuations}
\author{Ali \textsc{BOUFERROUM}}
\address[]{Ali BOUFERROUM, Laboratoire d'Analyse et de Math\'ematiques Appliqu\'ees, CNRS
 UMR8050, Universit\'e Paris-Est Marne-la-Vall\'ee, France.} 
\email{ali.bouferroum@univ-paris-est.fr} 
\date{ Compiled \today}
\subjclass[2010]{15B52; 60F05}
\begin{document}
\begin{center}
\maketitle 
\textbf{Universit\'e Paris-Est}
\end{center}
\begin{abstract}
 In this paper, we prove a universality result of convergence for a bivariate random process defined by the eigenvectors of a sample covariance matrix. Let $V_n=(v_{ij})_{i \leq n,\, j\leq m}$ be a $n\times m$ random matrix, where $(n/m)\to y > 0$ as $ n  \to \infty$, and let $X_n=(1/m) V_n V^{*}_n $ be the sample covariance matrix associated to $V_n \:$. Consider the spectral decomposition of $X_n$ given by $ U_n D_n U_n^{*}$, where $U_n=(u_{ij})_{n\times n}$ is an eigenmatrix of $X_n$. We prove, under some moments conditions, that the bivariate random process
$$
 \left( B_{s,t}^{n} =  \underset{1\leq j \leq \lfloor nt \rfloor}{\sum_{1\leq i \leq \lfloor ns \rfloor }} \left( |u_{i,j}|^2 - \frac{1}{n} \right) \right)_{(s,t)\in[0,1]^2}
$$
converges in distribution to a bivariate Brownian bridge. This type of result has been already proved for Wishart matrices (LOE/LUE) and Wigner matrices. This supports the intuition that the eigenmatrix of a sample covariance matrix is in a way "asymptotically Haar distributed". Our analysis follows closely the one of Benaych-Georges for Wigner matrices, itself inspired by Silverstein works on the eigenvectors of sample covariance matrices.

\vspace{0.5cm}
\noindent{\textsc{ Keywords and phrases.}} Random matrices; Sample covariance matrices; Haar measure; Eigenvectors; Delocalization; Brownian bridge; Spectral decomposition; Method of moments. \end{abstract} 

\section{Introduction}
The eigenvalues of random matrices attracted considerable attention in the recent years \cite{AGZ,BAISILV,BookTAO}. Less is known for the eigenvectors. Therefore, recent research on the limiting behavior of eigenvectors has attracted considerable interest among mathematicians and statisticians, see among others, Silverstein \cite{SILV,SIL1,SIL2}, Bai-Pan \cite{BaiPan}, Bai-Miao-Pan \cite{BaiMPan}, Ledoit-P\'ech\'e \cite{LedPec}, Benaych-Georges \cite{BENA}, Pillai-Yin \cite{PilYin}. The recent progress on the study of eigenvectors refers to a delocalization property shown for the eigenvectors of some types of random matrices, see Erd\"os-Schlein-Yau \cite{ErdYau,ErdSchYau}, Bordenave-Guionnet \cite{BorGui}, Schenker \cite{Schen}, Cacciapuoti-Maltsev-Schlein \cite{CAMASC}, Rudelson-Vershynin \cite{RUVE} and Vu-Wang \cite{VUWA}. For Wigner matrices, a universal properties of eigenvector coefficients were given recently, see Knowles-Yin \cite{KNOYIN} and Tao-Vu \cite{Tao_Vu}.\\

 In practical applications, the eigenvectors of large random matrices play a role as important as that played by the eigenvalues. For example, in multivariate analysis, the Principal Component Analysis is based on eigenvectors of  sample covariance matrices. The directions of the principal components are of particular interest,   
however, the exact distribution of the eigenmatrix (matrix of eigenvectors) of this type of matrices cannot be computed and few works had been devoted to this subject until quite recently. One of the reasons is that while the eigenvalues of an Hermitian matrix admit variational characterizations as extrema  of certain functions, the eigenvectors can be characterized as the argmax of these functions, hence are more sensitive to perturbations of the entries of the matrix.\\

Recently, it was proved in \cite{JIANG} that the entries of the first $ o(n/ \log n) $ columns of a Haar distributed matrix  can be approximated simultaneously by independent standard normals. Based on these evidentiary support and motivated by the fact that the eigenmatrix of Wishart matrix is Haar (uniformly) distributed, we believe that the eigenmatrix of a sample covariance matrix $X_n$ is "asymptotically Haar distributed" over the unitary group $\cU(n)$ of $n\times n$ unitary matrices for the complex case; or over the orthogonal group $\cO(n)$ of $n\times n$ orthogonal matrices for the real case. A question asked here is how to formulate the wording of "asymptotically Haar distributed"? Silverstein discussed this terminology in details in \cite{SILV}. \\

Let $V_n = (v_{ij})_{n \times m} \, , \: i= 1, \dots ,n \, ; \: j=1, \dots ,m=m(n)$ where $(n/m)\to y > 0$ as $ n  \to \infty$, be an observation matrix of i.i.d. real or complex random variables $\{v_{ij}\}_{i,j=1,2,\dots }$ such that 
 \begin{align}
 \E(v_{11}) = 0, \quad \quad  \E \left( |v_{11}|^{2} \right) =1, \hspace*{3cm}
\end{align} 
and $ V_j = (v_{1j} , \dots , v_{nj})^{'}$  be the $j^{th}$ column of  $V_n$. In this paper, we will consider a simplified version of sample covariance matrices with large dimension $n$ and sample size $m$
          $$X_n = (1/m) V_n V^{*}_n, $$ 
 where $ V^{*}_n $ denotes the conjugate transpose of the matrix $V_n$. Let us define the cumulative distribution of the eigenvalues of $X_n$, for each $u \in \R$, as 
  $$  F^X_n(u)= \frac{1}{n}\sum_{i=1}^{n} \textbf{1}_{\lbrace\gl_i\leq u\rbrace}, $$ 
this function describes the global behavior of the spectrum of $X_{n}$. Recall that for a matrix $X_n$ defined as above, the previous empirical cumulative function $F^X_n$ converges almost surely for every $u \geq 0$, as $ n  \to \infty$, to a non-random distribution function $F^{MP}_y$ which has the Marchenko-Pastur density 
$$ f^{MP}_y(u)= \left(1-\frac{1}{y}\right)_{+} \gd_{0} + \frac{\sqrt{\left(b - u \right)\left( u-a \right)}}{2\pi y u} \: \mathbf{1}_{[a,b]}(u), $$
where $ a = (1-\sqrt{y})^2$ and $ b  = (1+\sqrt{y})^2 $ (atom $1- 1/y$ at the origin if and only if $y>1$), see \cite{YIN} and \cite[Theorems 3.6, 3.7]{BAISILV}\label{MP}.\\

Let $ U_n D_n U_n^{*}$ denote the spectral decomposition of the sample covariance matrix $X_n$, where $D_n = \text{diag}(\gl_1, \dots , \gl_n)$ , and the $ \gl_i$'s   are the eigenvalues of $X_n$ arranged along the diagonal of $D_n$ in non-decreasing order, and  $U_n=\{u_{ij}\}$ is the associated eigenmatrix for $X_n$. Let us define a bivariate random process $B_{s,t}^{n}$ by 

\begin{align}
\left( B_{s,t}^{n} = \sqrt{\frac{\gb}{2}} \underset{1\leq j \leq \lfloor nt \rfloor}{\sum_{1\leq i \leq \lfloor ns \rfloor }} \left( |u_{i,j}|^2 - \frac{1}{n} \right) \right)_{(s,t)\in[0,1]^2}
\end{align}
where  $\gb = 2 $ in the complex case and $\gb = 1 $  in the real case, and $ \lfloor a \rfloor$ denotes the greatest integer less than or equal to a. It is well known \cite{Don_Rou} that if $U_n$ is Haar distributed over the group $\cU(n)$  or the group $\cO(n)$, then $B_{s,t}^{n}$ weakly converges (i.e. in the sense of convergence of all finite-dimensional marginals) to a Brownian bridge $(B_{s,t})$ as $n$ tends to infinity, i.e: the centered continuous Gaussian
process $(B_{s,t})_{(s,t)\in [0,1]^2}$ with covariance 
\begin{align}
\E\left(B_{s,t} B_{s^{'} , t^{'}} \right)= \left( \min \{ s, s^{'} \}- s s^{'} \right) \left( \min \{ t, t^{'} \}- t t^{'} \right).
\end{align}
 Conversely, if $B_{s,t}^{n}$ weakly converges to a Brownian bridge, it then reveals some evidence supporting the conjecture that the eigenmatrix $U_n$ is asymptotically Haar distributed.\\
 
 In this paper, we will prove that for a sample covariance matrix $X_n$ defined as above, $B_{s,t}^{n}$ has a limit in a weaker sense if $v_{11}$ has moments of all orders, and that this weak limit is the bivariate Brownian bridge if and only if $v_{11}$ has the same fourth moment as in the case of LOE/LUE matrix (Wishart-Laguerre orthogonal/unitary ensembles).  This work is inspired by the work of Benaych-Georges for Wigner matrices \cite{BENA}, itself is inspired by Silverstein's works in \cite{SIL1,SIL2} for a univariate process defined by the eigenmatrix of a sample covariance matrix as
 $$ Y^{n}_t = \sqrt{\frac{n}{2}} \sum_{1\leq i \leq \lfloor nt \rfloor } \left( |y_{i}|^2 - \frac{1}{n} \right),  \quad \quad \mathrm{where \quad} y=X_nx_n \mathrm{\quad for \: some \: vector \: x_n }. $$

The rest of the paper is organized as follows. The main theorem is presented in section \ref{sec2} with some remarks.  The proof of this theorem is mainly contained in Sections \ref{secoutline} and \ref{secproof}: Essentially, the problem is transformed into showing convergence on an appropriate space $D[0,1]\times D[0,+\infty[$ instead of the space $D[0,1]^2$. After that, the proof will consist of studying the moments of a weighted spectral law of $X_n$ according to the process $B_n$. In section \ref{tightness}, we finish by giving a version of tightness and convergence in the Skorokhod topology of the process $B^n_{s,t}$ under some additional hypotheses on the atom distribution.  \\
\vspace*{0.4cm}

\noindent{\bf Acknowledgment: } I would like to express my sincere appreciation and gratitude to my advisor Djalil Chafa\"i for his academic guidance and enthusiastic encouragement. I also would like to thank Florent Benaych-Georges for pointing out some references and for his encouragement.

\section{Main result}  \label{sec2}
Let us consider a matrix  $V_n = (v_{ij}^{(n)})_{n \times m} \, , \: i= 1, \dots ,n \, ; \: j=1, \dots ,m=m(n)$ where $(n/m)\to y > 0$ as $ n $ tends to infinity. Let $X_n = (1/m) V_n V^{*}_n$ be its associated sample covariance matrix of dimension $n$ and sample size $m$. Let $ U_n D_n U_n^{*}$ denote the spectral decomposition of the sample covariance matrix $X_n$, where $D_n = \text{diag}(\gl_1, \dots , \gl_n)$ , and $ \gl_i$'s are the eigenvalues of $X_n$ arranged along the diagonal of $D_n$ with a non-decreasing order, and  $U_n=\{u_{ij}\}_{1 \leq i,j\leq n}$ is the associated eigenmatrix of $X_n$.
Note that $U_n$ is not uniquely defined, however, one can choose it in any measurable way. We consider the bivariate c\`{a}d-l\`{a}g process $\left(B_{s,t}^{n}\right)_{{(s,t)\in[0,1]^2}}$ defined as:
\[
 \left( B_{s,t}^{n} = \sqrt{\frac{\gb}{2}} \underset{1\leq j \leq \lfloor nt \rfloor}{\sum_{1\leq i \leq \lfloor ns \rfloor }} \left( |u_{i,j}|^2 - \frac{1}{n} \right) \right)_{(s,t)\in[0,1]^2},
\]
where  $\gb = 2 $ in the the complex case and $\gb = 1 $  in the real case.

\vspace{0.3cm}
\begin{theorem}[Main result] \label{th1}
Suppose in the definition above of the sample covariance matrix $X_n= (1/m) V_n V^{*}_n $ that 
\begin{align} \label{Ass0}
\{v_{ij}^{(n)}\}_{i,j=1,2,\dots} \quad are \:\, i.i.d. \:\, complex \:\, or \:\, real \:\, random \:\, variables,
\end{align}
with
\begin{align} \label{Ass1}
 \E(v_{11}^{(n)}) = 0, \quad \quad  \E |v_{11}^{(n)}|^{2} =1, 
 \end{align}
and  
 \begin{align} \label{Ass2}
 \forall \, k \geq 0, \: \; \sup_n \, \E|v_{1,1}^{(n)}|^{k}< \infty ,
 \end{align}
Then the sequence 
$$ \left(\dis(B^{n})\right)_{n\geq 1} $$ 
has a unique possible accumulation point supported by $C(\left[0,1\right]^2)$ in the sense of convergence of all finite-dimensional marginals. This accumulation point is the distribution of a centered Gaussian process which depends on the distribution of $v_{11}^{(n)}$ only through $\lim_{ n  \to \infty}  \E \left( |v_{11}^{(n)}|^{4} \right)$, and which is the bivariate Brownian bridge if and only if
 \begin{align}
 \lim_{ n  \to \infty}  \E \left( |v_{11}^{(n)}|^{4}\right)=4-\gb.
 \end{align} 
\end{theorem}
\vspace{0.3cm}
\begin{rem}[Dependence of entries on $n$]
The distribution of the entries $v_{i,j}=v_{i,j}^{(n)}$ are allowed to depend on $n$. For brevity of notations, we write $v_{i,j}$ instead of $v_{i,j}^{(n)}$.  
\end{rem}
\begin{rem}[Matching with LOE/LUE] Note that the unique possible accumulation point supported by $C([0,1]^2) $ of our sequence $(\dis(B^{n}))_{n\geq 1}$ which is a centered Gaussian process depends on the distribution of the $v_{i,j}$'s only through $ \, \lim_{ n  \to \infty}  \E ( |v_{11}|^{4} )$, and this limiting distribution is the bivariate Brownian bridge if and only if $ \, \lim_{ n  \to \infty} \: \E ( |v_{11}|^{4} )\: $ is the same as for a LOE or LUE matrix, i.e. equal to $4-\gb$. 
\end{rem}

\newpage
 
 \section*{Simulation}
    \begin{figure}[!h]
\begin{center}
\includegraphics[width=15cm,height=7cm]{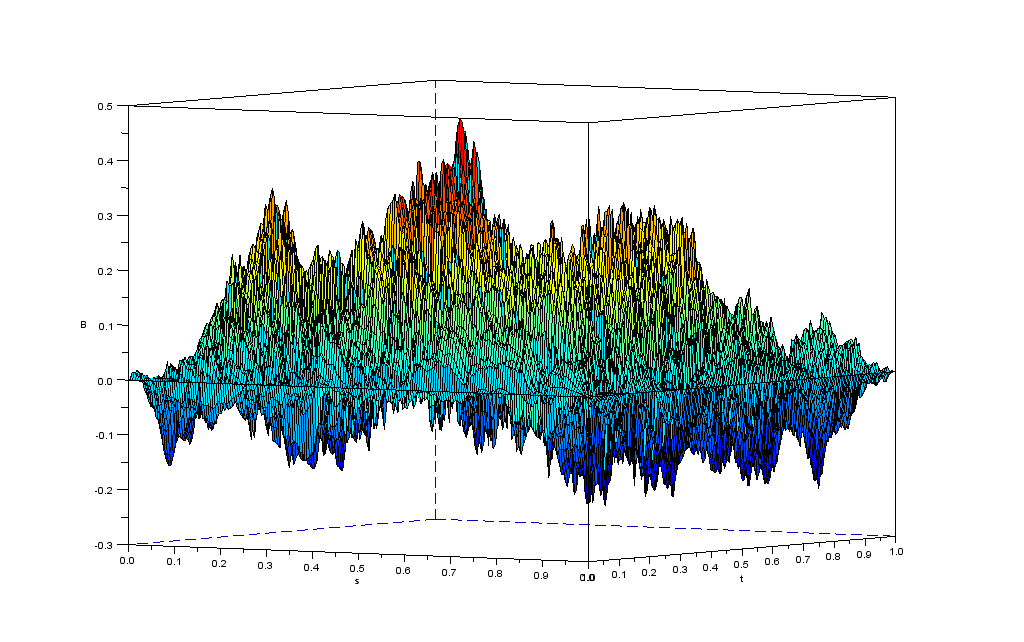}
\caption{ Simulation of the bivariate random process $B^n(s,t)$ for a sample covariance matrix with standard real normal atom distribution (Wishart marix). the matrix $V_n$ is $n \times m$ with $n=m=500$.}  
\end{center}
\end{figure}

\begin{figure}[!h]
\begin{center}
\includegraphics[width=15cm,height=7cm]{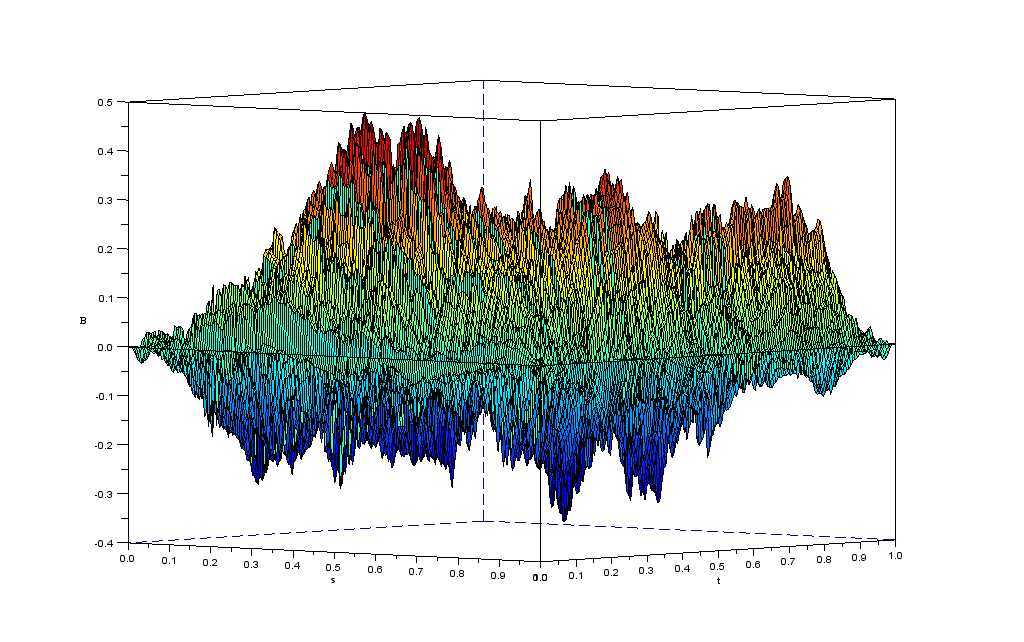}
\caption{ Simulation of the bivariate random process $B^n(s,t)$ for a sample covariance matrix with atom distribution: $v_{1,1}=CB-\E(CB)$ with $CB \stackrel{d}{=}\text{Binomial}(6,(1/2-\sqrt{1/12}))$. The matrix $V_n$ is $n \times m$ with $n=m=500$.} 
\end{center}
\end{figure}

\begin{figure}[!h]
\begin{center}
\includegraphics[width=15cm,height=7cm]{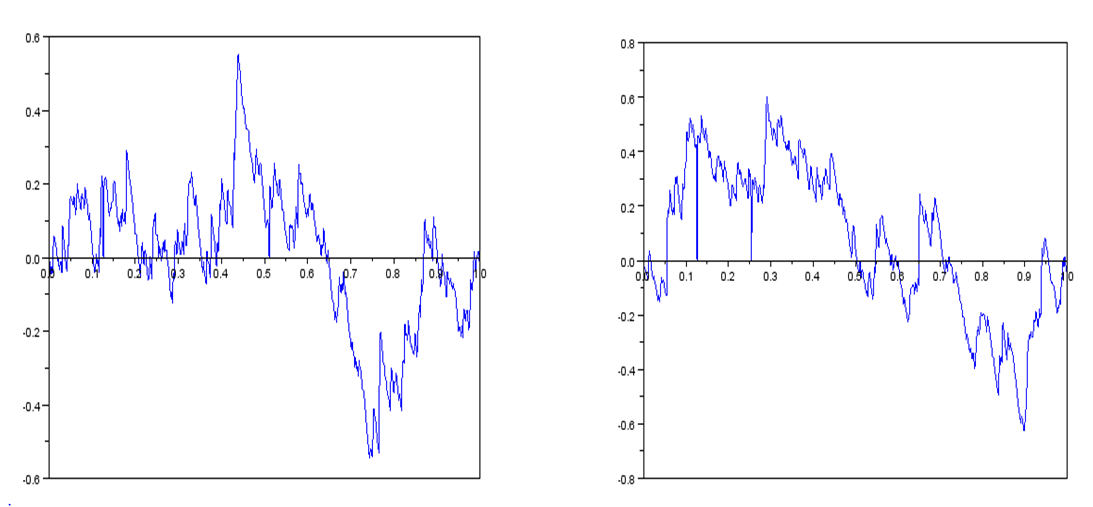}
\caption{ Simulation of the univariate random process $B^n(s,1)$ for two different choices of  atom distributions: 
 {\bf Left:} Wishart matrix.  {\bf Right:}  Sample covariance matrix with atom distribution $v_{1,1}=CB-\E(CB)$ with $CB\stackrel{d}{=}\text{Binomial}(6,(1/2-\sqrt{1/12}))$.  For both pictures, the matrix $V_n$ is $n\times m$ with $n=m=500$.}
\end{center}
\end{figure}

\newpage

\section{Outline of the proof of Theorem \ref{th1}}  \label{secoutline}
In this paper we denote:
\begin{enumerate}
\item[$\bullet$]  $C([0,1]^2)$ $\left(\text{resp}. \quad C_c([0,1] \times [0, + \infty) \right) $ the space of real valued continuous functions on $[0,1]^2$ (resp. of real valued compactly supported continuous functions on $[0,1] \times [0, +\infty[$), endowed with the uniform convergence topology.\\
$\,$
\item[$\bullet$] $D_c(\R,[0,1])$ the set of compactly supported c\`{a}d-l\`{a}g functions on $\R$ taking values in $[0,1]$, endowed with the topology defined by the fact that $ f_n \to f $ if and only if the
bounds of the support of $f_n$ tend to those of the support of $f$ and for all $M > 0$, after restriction to $[-M,M]$, $ f_n \to f $ with the topology of $D[-M,M]$ being deduced from the Skorokhod topology of $D([0, 1])$ defined in \cite[Chapter 3]{BILL}.\\
$\,$
\item[$\bullet$]  $D([0,1]^2)$ $\left(\text{resp}. \quad D_c([0,1] \times [0, +\infty[ \right) $ the space of real valued functions $f:[0,1]^2 \to \R $ (resp. of compactly supported functions $f:[0,1] \times [0, +\infty[ \to \R$) admitting limits in all "orthants", more precisely such that for each $s_0 , t_0$, for each pair of symbols $\circ ,  \bullet \in \lbrace <, \geq \rbrace $,
$$ \underset{t \to_{\bullet}\, t_0 }{\underset{s \to_{\circ}\, s_0}{\lim}}f(s,t) $$
exists, and is equal to $f(s_0,t_0)$ if both $\, \circ$ and $\bullet \, $ are $\geq$. The space $D([0,1]^2)$ is endowed with the Skorokhod topology defined in \cite{BIWI} and the space $D_c([0,1] \times [0, + \infty[)$ is endowed with the topology defined by: $f_n\to f$ if and only if for all $M > 0$, after restriction to $[0,1] \times [0,M]$, $f_n\to f$ in the sense of the space $D([0,1]^2)$.\\
$\,$
\item[$\bullet$] $D_0([0,1]^2)$ the set of functions in $D([0,1]^2)$ vanishing at the border of $[0, 1]^2$, endowed with the induced topology.
 \end{enumerate}
 \subsection{From $D([0,1]^2)$ to $D([0,1] \times [0, +\infty[ )$} 
As we have seen in the introduction, the cumulative distribution function $F^X_n$ of $X_n$ converges almost surely, as $ n  \to \infty$, to a non-random distribution function $F^{MP}_y$ (Marchenko-Pastur law). The proof of Theorem \ref{th1} can be reduced to the following remark, inspired by some ideas of Silverstein \cite{SIL1,SIL2} and of Benaych-Georges \cite{BENA}: even though we do not have any "direct access" to the eigenvectors of $X_n$, we have access to the process  $ \left( B^{n}_{s,F^X_n(u)}\right)_{s\in[0,1],u \in  [0, +\infty[}$, for $ F^X_n(u)= \frac{1}{n}\sum_{i=1}^{n} \textbf{1}_{\lbrace\gl_i\leq u\rbrace}$. Indeed,  
	   $$ B^{n}_{s,F^X_n(u)} = \sqrt{\frac{\gb}{2}} \sum_{1\leq i \leq \lfloor ns \rfloor } \sum_{1\leq j \leq n:\gl_j \leq u } \left( |u_{i,j}|^2 - \frac{1}{n} \right),$$
hence, for all fixed $ s \in [0, 1]$ , the function $u \in \R \mapsto B^{n}_{s,F^X_n(u)}$ is the cumulative distribution function of the signed measure 
\begin{align} \label{LSP}
\sqrt{\frac{\gb}{2}} \sum_{1\leq i \leq \lfloor ns \rfloor } \sum_{1\leq j \leq n} \left( |u_{i,j}|^2 - \frac{1}{n} \right)\gd_{\gl_j},
\end{align}
which can be considered as a difference between two random  probability measures:\\
 $ \sum_{1\leq j \leq n}  |u_{i,j}|^2 \, \gd_{\gl_j}$ (weighted spectral law of $X_n$) and  $ \frac{1}{n} \sum_{1\leq j \leq n} \gd_{\gl_j}$ (empirical spectral law of $X_n$). The law (\ref{LSP}) can be studied via its moments  
$$ \sum_{1\leq i \leq \lfloor ns \rfloor } \left(e_i^{*}X_n^k e_i - \frac{1}{n} \text{Tr} X_n^k \right) $$
for $k \geq 1$, the $e_i$'s being the vectors of the canonical basis. From the asymptotic behavior of the moments of the signed measure (\ref{LSP}), one can then find out the asymptotic behavior of its cumulative distribution function.\\
Once the asymptotic distribution of the process $ \left( B^{n}_{s,F^X_n(u)}\right)_{s\in[0,1],u \in  [0, +\infty[}$ is identified, one can obtain the asymptotic distribution of the process $(B^n_{ s,t})_{(s,t)\in [0,1]^2}$, because $$F^X_n(u) \underset{n \to \infty}{\longrightarrow} F^{MP}_y(u), \quad \textrm{almost surely for every} \: u \geq 0.$$ 
 
The following proposition is the key of the proof, since it allows transferring our problem from the eigenvectors to some more accessible objects: the weighted spectral distributions of the sample covariance matrix $X_n$.

\begin{prop}[From the process $B_n$ to a weighted spectral process] \label{prop1} To prove Theorem \ref{th1}, it suffices to prove that each finite-dimensional marginal distribution of the process 
$$ \left( \sum_{1\leq i \leq \lfloor ns \rfloor } \left(e_i^{*}X_n^k e_i - \frac{1}{n} \text{Tr} X_n^k \right)\right)_{s\in [0,1], k \geq 1 } $$
converges to a centered Gaussian process and that the covariance of the limiting process depends
on the distribution of the $v_{i,j}$'s only through  $ \, \lim_{ n  \to \infty}  \E ( |v_{11}|^{4} )$, and that this covariance is the one of the bivariate Brownian bridge if and only if $\: \lim_{ n  \to \infty} \E ( |v_{11}|^{4} )= 4-\gb $.
\end{prop}

\begin{Proof} It is known \cite{YIN,BAISILV} that the cumulative distribution function $F^X_n(u)= \frac{1}{n}\sum_{i=1}^{n} \textbf{1}_{\lbrace\gl_i\leq u\rbrace} $ of the matrix $X_n$ converges almost surely,  
as $n$ tends to infinity, to a non-random distribution function $F^{MP}_y$ defined by means of the Marchenko-Pastur law. Since the limit is continuous and compactly supported on $(0,+\infty)$, this convergence is uniform
$$  \sup_{u\in[0,+\infty[} |F^X_n(u)-F^{MP}_y(u)| \underset{n \to \infty}{\longrightarrow} 0    \quad \quad \text{almost surely}. $$
Hence, it follows that 
$$F^X_n  \underset{n \to \infty}{\longrightarrow}  F^{MP}_y   \quad \quad \mathrm{in} \:\: D_c([0, +\infty),[0,1]), $$
see \cite[section 10.1.2]{BAISILV}. Moreover, the map 
   $$ D_0([0,1]^2) \times D_c([0, +\infty),[0,1]) \to D_c([0,1] \times [0, +\infty)) $$
   $$ ((G_{s,t})_{s,t \in [0,1]}, (g(u))_{u\in [0, +\infty)}) \mapsto (G_{s,g(u)})_{(s,u)\in [0,1] \times [0, +\infty)}$$ 
is continuous at any pair of continuous functions. Hence for any continuous process $(B_{s,t})_{s,t \in [0,1]}$ whose distribution is an accumulation point of the sequence	$(\dis(B^{n}))_{n\geq 1}$ for the Skorokhod
topology in $D([0, 1]^2)$, the process
 $$ \left(B^n_{s,F^{X}_n (u)} \right)_{s,u \in [0,1] \times [0, +\infty)} $$
 converges in distribution (up to the extraction of a subsequence) to the process
 $$ \left(B_{s,F^{MP}_y(u)} \right)_{s,u \in [0,1] \times [0, +\infty)}. $$
This assertion relies on two results which can be found in \cite[Theorem 4.4 and Corollary 1 of Theorem 5.1 in Chapter 1]{BILL}. Now, note that $F^{MP}_y:[0, +\infty) \to [0,1]$ admits a right inverse, so the distribution of the process $(B_{s,t})_{s,t} \in [0,1]$ is entirely determined by that of the process  $(B_{s,F^{MP}_y(u)})_{s,u \in [0,1] \times [0, +\infty)} $. Therefore, to prove Theorem \ref{th1}, it suffices to prove that the sequence
\begin{align} \label{fordis}
\left( \dis (B^n_{s,F^{X}_n (u)})_{s,u \in [0,1] \times [0, +\infty)})\right)_{n \geq 1}
\end{align}
converges to a centered Gaussian process which depends on the distribution of the $v_{i,j}$'s only through $\: \lim_{ n  \to \infty}  \E \left( |v_{11}|^{4} \right)$, and which is the bivariate Brownian bridge if and only if \\
$$ \lim_{ n  \to \infty}  \E \left( |v_{11}|^{4}\right)=4-\gb.$$

Now, let us prove that any (random) function  $f \in C_c([0,1] \times [0, +\infty))$ is entirely determined by the collection of real numbers $(\int_{u\in [0, +\infty)} u^k f(s,u)  \, du )_{s \in [0,1], k \geq 0}$.
\begin{lemma}[Technical characterization] Let $f$ be a random variable in $C_c([0,1] \times [0, +\infty))$ such that with probability
one, $f(s, u) = 0$ when $u > r, r > 0$. Then the distribution of f is entirely determined by the finite-dimensional marginals of the process  
\begin{align} \label{lemgau}
\left(\int_{u\in [0, +\infty)} u^k f(s,u)  \, du \right)_{s \in [0,1], k \geq 0}.
\end{align} 
Moreover, in the case where the finite-dimensional marginals of the process of (\ref{lemgau}) are Gaussian
and centered, then so are those of $f$.
\end{lemma}

\begin{Proof} Let us fix $(s, u_0) \in [0,1] \times [0,r]$ and let, for each $p \geq 1$, $(P_{p,q})_{q \geq 1}$ be a sequence of polynomials that is uniformly bounded on $[0, r+1]$ and that converges pointwise to  $1_{[u_0,u_0 +1/p]}$
on $[0; r+1]$. Then one has, with probability one,
$$f(s,u_0)= \lim_{p\to \infty } p \int_{u_0}^{u_0 +1/p} f(s,u) du = \lim_{p\to \infty }\lim_{q\to \infty } p \int_{u\in [0, +\infty)} P_{p,q} (u) f(s,u) du. $$
This proves the lemma, because any almost sure limit of a sequence of variables belonging to a
space of centered Gaussian variables is Gaussian and centered.
\end{Proof}
Since the fourth moment of the entries of $V_n$ is finite, we know that the largest eigenvalue of the  sample covariance matrix $X_n$ converges, almost surely, to $ b  = (1+\sqrt{y})^2 $ (see \cite{BSY,BY}). Hence, for any random variable $f$ taking values in $C_c([0, 1] \times [0,+\infty))$ such that the distribution of $f$ is a limiting point of the sequence of (\ref{fordis}), we have $f(s, u) = 0$, almost surely,  when $u> b+\epsilon$.\\
As a consequence, it follows from the previous lemma and from what precedes that in order to prove Theorem \ref{th1}, it suffices to prove that each finite-dimensional marginal distribution of the process	 
$$ \int_{u \in \R} u^{k} B^{n}_{s,F^X_n(u)} \, du $$ 
converges to a centered Gaussian measure and that the covariance of the limit process depends on the distribution of the $v_{i,j}$'s s only through $ \, \lim_{ n  \to \infty}  \E ( |v_{11}|^{4} )$, and that this covariance is the one of the bivariate Brownian bridge if and only if $\: \lim_{ n  \to \infty} \E ( |v_{11}|^{4} )= 4-\gb $.\\

 Recall that :
 \[
 \begin{split} 
 B^{n}_{s,F^X_n(u)} & = \sqrt{\frac{\gb}{2}} \sum_{1\leq i \leq \lfloor ns \rfloor } \sum_{1\leq j \leq n:\gl_j \leq u } \left( |u_{i,j}|^2 - \frac{1}{n} \right)\\
                    & = \sqrt{\frac{\gb}{2}} \sum_{1\leq i \leq \lfloor ns \rfloor } F_{\mu_{X_n,e_i}-\mu_{X_n}},
 \end{split}
\]
where
\begin{enumerate}
\item[$\bullet$] $\mu_{X_n}$ is the empirical spectral law of $X_n$. 
\item[$\bullet$] $\mu_{X_n,e_i}$ is the weighted spectral law of $X_n$, defined by $\mu_{X_n,e_i}= \sum_{j=1}^{n} |u_{i,j}|^2 \delta_{\gl_j}$.
\item[$\bullet$] $F_{\mu_{X_n,e_i}-\mu_{X_n}} $ is the cumulative distribution function of the null-mass signed measure $\mu_{X_n,e_i}-\mu_{X_n}$.
\end{enumerate}
Now, let us give the following lemma to complete the proof of Proposition \ref{prop1}.
\begin{lemma}[Moment's calculation rule]  
Let $\mu$ be a compactly supported null-mass signed measure and set $F_{\mu}(u) =
\mu((- \infty, u])$. Then for all $k \geq 0$,
 $$ \int_{u\in \R} u^k  F_{\mu}(u) du = - \int_{x\in \R} \frac{x^{k+1}}{k+1} d\mu(x). $$
\end{lemma}
\begin{Proof}
Let $a<b$ be such that the support of $\mu$ is contained in the open interval $(a, b)$. $F_{\mu}$ is
null outside $(a, b)$ and satisfies $F_{\mu}(u) = - \mu((u, b))$, so by Fubini's Theorem,
\[
\begin{split}
\int_{u\in \R}  u^k F_\mu(u) \, du & =   \int_{x=a}^b \int_{u=a}^{x} -u^k \, du \, d\mu(x) \\
& =\int_{x=a}^b \frac{-x^{k+1}}{k+1} \, d\mu(x)+ \frac{a^{k+1}}{k+1}\mu((a,b))\\
& =\int_{x=a}^b \frac{-x^{k+1}}{k+1}d\mu(x).
\end{split}
\]
\end{Proof}
 It follows from this lemma  that for all $s\in [0,1]$, $k\geq 0$,
  \[
  \begin{split}
  \int_{u\in \R} u^k B^n_{s,F^X_n}(u) \, du & = - \int_{x\in \R}  \frac{-x^{k+1}}{k+1} \,  d(\mu_{X_n,e_i}-\mu_{X_n})(x)\\
& = \frac{-1}{k+1}\sqrt{\frac{\beta}{2}} \sum_{1\leq i \leq \lfloor ns \rfloor } \sum_{1\leq j \leq n} \left( |u_{i,j}|^2 - \frac{1}{n} \right) \gl_j^{k+1}\\
  & = \frac{-1}{k+1}\sqrt{\frac{\beta}{2}}\sum_{1\leq i\leq \lfloor ns \rfloor }\left(e_i^*X_n^{k+1}e_i-\frac{1}{n} \text{Tr} X_n^{k+1}\right), 
\end{split}
\]
which completes the proof of Proposition \ref{prop1}.
\end{Proof}

It follows from all what precedes that Theorem \ref{th1} is a direct consequence of the following proposition, whose proof is in Section \ref{secproof}. 
\begin{prop}[Convergence of Moments]\label{mainprop}
Under Assumptions (\ref{Ass0}), (\ref{Ass1}) and (\ref{Ass2}) in Theorem \ref{th1}, each finite-dimensional marginal distribution
of the process
 $$ \left( \sum_{1\leq i \leq \lfloor ns \rfloor } \left(e_i^{*}X_n^k e_i - \frac{1}{n} \text{Tr} X_n^k \right)\right)_{s\in [0,1], k \geq 1 } $$
 converges to   a centered  Gaussian measure. The covariance of the limit distribution,  denoted by  $$\left(\Cov_{s_1,s_2}(k_1,k_2)\right)_{s_1,s_2\in [0,1], k_1,k_2 \ge 1},$$  depends on the distribution of the $v_{i,j}$'s only through  $\lim_{n \to \infty}\E[|v_{1,1}|^4]$. Moreover, we have
 $$ \Cov_{s_1,s_2}(1,1)=(\lim_{n\to \infty} \E[|x_{1,2}|^4]-1)(\min\{s_1,s_2\}-s_1s_2).$$
\end{prop} 
\section{Proof of Proposition \ref{mainprop}} \label{secproof}
Note that the expectation of the weighted spectral law $\mu_{X_n,e_i}=\sum_{j=1}^{n} |u_{i,j}|^2 \delta_{\gl_j}$ does not depend on $i$. So for all $s \in [0,1]$, $k \geq  1$,
\begin{equation} \label{centrage}
\sum_{1\leq i\leq ns}(e^*_iX_n^ke_i-\frac{1}{n} \text{Tr}(X_n^k))=\sum_{1\le i\le ns}(e^*_iX_n^ke_i-\E[e^*_iX_n^ke_i])-\frac{\lfloor ns\rfloor}{n}\sum_{1\leq i\leq n}(e^*_iX_n^ke_i-\E[e^*_iX_n^ke_i]).
\end{equation}
Therefore, we are led to study the limit, as $n \to \infty$, of the finite-dimensional marginal distributions of the process
\begin{align} \label{mainfor}
\left(\sum_{1\leq i\leq ns}(e^*_iX_n^ke_i-\E[e^*_iX_n^ke_i])\right)_{s\in [0,1], k \geq  1}.
\end{align}
Let us fix $p \geq 1$, $s_1, \ldots, s_p\in [0,1]$ and $k_1, \ldots, k_p \geq 1$. We shall study the limit, as $n$ tends to infinity, of 
 \begin{align} \label{EEE}
 \E \left[ \prod_{l=1}^p \sum_{1\leq i\leq ns_l}(e^*_iX_n^{k_l} e_i-\E[e^*_iX_n^{k_l}e_i]) \right]
 \end{align}
 We introduce the set
\begin{align} \label{EE}
\cE:=\{\ovl{0}^1, \ldots,\ovl{k_1}^1\} \cup  \cdots\cdots \cup \{\ovl{0}^p, \ldots,\ovl{k_p}^p\},
\end{align} 
where the sets $\{\ovl{0}^1, \ovl{1}^1, \ldots\}$, $\{\ovl{0}^2, \ovl{1}^2, \ldots\}$, \ldots,  $\{\ovl{0}^p, \ovl{1}^p, \ldots\}$  are $p$ disjoint copies of the set of non-negative integers. The set $\cE$ is ordered as presented in (\ref{EE}). In the rest of this paper, we denote $(1,\ldots, n)^k$ the set of k$-$tuples of a set $\{ 1,2, \ldots, n \} $.\\
The expectation (\ref{EEE}) can be expanded and expressed as a sum on the set $(1,\ldots, n)^{K}$ indexed by the set $\cE$ introduced above, where  $K= k_1+ \cdots + k_p$ . We get 
\begin{align} \label{ME}
 \E \left[ \prod_{l=1}^p \sum_{1\leq i\leq ns_l}(e^*_iX_n^{k_l} e_i-\E[e^*_iX_n^{k_l}e_i]) \right] 
 \end{align}
 $$
= \sum_{\pi\in (1,\ldots,n)^{K}}\E \left[ \prod_{l=1}^p \left(x_{\pi(\ovl{0}^l),\pi(\ovl{1}^l)}\cdots x_{\pi(\ovl{k_l -1}^l),\pi(\ovl{k_l}^l)} - \E(x_{\pi(\ovl{0}^l),\pi(\ovl{1}^l)}\cdots x_{\pi(\ovl{k_l -1}^l),\pi(\ovl{k_l}^l)})\right) \right] \mathbf{1}_{\{1\leq \pi(\ovl{0}^l) = \pi(\ovl{k_l}^l) \leq ns_l \}}, 
$$
As we have not sufficient information on the laws of the $ x_{i,j}$'s , we need to write the previous expression in terms of elements of the matrix $V_n$. Let us consider
$$ \text{Pr} = x_{\pi(\ovl{0}^l),\pi(\ovl{1}^l)} x_{\pi(\ovl{1}^l),\pi(\ovl{2}^l)} \cdots x_{\pi(\ovl{k_l -1}^l),\pi(\ovl{k_l}^l)}.$$

We have
\[
\begin{split}
\text{Pr} & = \frac{1}{m^{k_l}} \left[\sum_{j=1}^m v_{\pi(\ovl{0}^l),j} v^*_{j,\pi(\ovl{1}^l)}\right]\left[\sum_{j=1}^m v_{\pi(\ovl{1}^l),j} v^*_{j,\pi(\ovl{2}^l)}\right] \cdots \left[\sum_{j=1}^m v_{\pi(\ovl{k_l -1}^l),j} v^*_{j,\ovl{k_l}^l)}\right]\\
 & = \frac{1}{m^{k_l}} \sum_{\gga\in(1,\ldots,m)^{k_l}} v_{\pi(\ovl{0}^l),\gga(\ovl{1}^l)}v^*_{\gga(\ovl{1}^l),\pi(\ovl{1}^l)}v_{\pi(\ovl{1}^l),\gga(\ovl{2}^l)}v^*_{\gga(\ovl{2}^l),\pi(\ovl{2}^l)}\cdots v_{\pi(\ovl{k_l -1}^l),\gga(\ovl{k_l}^l)}v^*_{\gga(\ovl{k_l}^l)\pi(\ovl{k_l}^l)}. 
\end{split}
\]
Thus, we obtain that
 \[
 \E \left[ \prod_{l=1}^p \sum_{1\leq i\leq ns_l}(e^*_iX_n^{k_l} e_i-\E[e^*_iX_n^{k_l}e_i]) \right]  \]
\begin{equation} \label{MME}
= m^{-K} \sum_{\pi\in (1,\ldots,n)^{K}}\E \left[ \prod_{l=1}^p  \sum_{\gga\in(1,\ldots,m)^{k_l}} \left( V^{k_l}_{l,\pi,\gga}- \E \, V^{k_l}_{l,\pi,\gga} \right)\right]\mathbf{1}_{\{1\leq \pi(\ovl{0}^l) = \pi(\ovl{k_l}^l) \leq ns_l \}},
\end{equation} 
where
$$V^{k_l}_{l,\pi,\gga} = v_{\pi(\ovl{0}^l),\gga(\ovl{1}^l)}v^*_{\gga(\ovl{1}^l),\pi(\ovl{1}^l)}v_{\pi(\ovl{1}^l),\gga(\ovl{2}^l)}v^*_{\gga(\ovl{2}^l),\pi(\ovl{2}^l)}\cdots v_{\pi(\ovl{k_l -1}^l),\gga(\ovl{k_l}^l)}v^*_{\gga(\ovl{k_l}^l)\pi(\ovl{k_l}^l)}.
$$
The product and the sum in the expectation of (\ref{MME}) can be expressed and developed as the following: 
$$ \prod_{l=1}^p  \sum_{\gga\in(1,\ldots,m)^{k_l}} v_{\pi(\ovl{0}^l),\gga(\ovl{1}^l)}v^*_{\gga(\ovl{1}^l),\pi(\ovl{1}^l)}v_{\pi(\ovl{1}^l),\gga(\ovl{2}^l)}v^*_{\gga(\ovl{2}^l),\pi(\ovl{2}^l)}\cdots v_{\pi(\ovl{k_l -1}^l),\gga(\ovl{k_l}^l)}v^*_{\gga(\ovl{k_l}^l)\pi(\ovl{k_l}^l)} $$
$$= \sum_{\gga\in(1,\ldots,m)^{K}} (v_{\pi(\ovl{0}^1),\gga(\ovl{1}^1)}v^*_{\gga(\ovl{1}^1),\pi(\ovl{1}^1)} \cdots  v_{(\ovl{k_1 -1}^1),\gga(\ovl{k_1}^1)}v^*_{\gga(\ovl{k_1}^1)\pi(\ovl{k_1}^1)}) \cdots  (v_{\pi(\ovl{0}^p),\gga(\ovl{1}^p)}v^*_{\gga(\ovl{1}^p),\pi(\ovl{1}^p)} \cdots  v^*_{\gga(\ovl{k_p}^p)\pi(\ovl{k_p}^p)}) $$
$$= \sum_{\gga\in(1,\ldots,m)^{K}} \prod_{l=1}^p \left(v_{\pi(\ovl{0}^l),\gga(\ovl{1}^l)}v^*_{\gga(\ovl{1}^l),\pi(\ovl{1}^l)} \cdots  v_{\pi(\ovl{k_l -1}^l),\gga(\ovl{k_l}^l)}v^*_{\gga(\ovl{k_l}^l)\pi(\ovl{k_l}^l)}\right).
$$
Therefore, we find that the quantity (\ref{ME}) is equal to:
\begin{align} \label{arrang}
m^{-K} \sum_{\pi\in (1,\ldots,n)^{K}}\sum_{\gga\in(1,\ldots,m)^{K}} \E \left[\prod_{l=1}^p \left(V_{l,\pi,\gga}-\E \, V_{l,\pi,\gga}\right)\right]\mathbf{1}_{\{1\leq \pi(\ovl{0}^l) = \pi(\ovl{k_l}^l) \leq ns_l \}},
\end{align}  
where 
$$ V_{l,\pi,\gga} = v_{\pi(\ovl{0}^l),\gga(\ovl{1}^l)}v^*_{\gga(\ovl{1}^l),\pi(\ovl{1}^l)} \cdots  v_{\pi(\ovl{k_l -1}^l),\gga(\ovl{k_l}^l)}v^*_{\gga(\ovl{k_l}^l)\pi(\ovl{k_l}^l)}. 
$$
$\,$\\

Now, as we work directly with the elements $v_{i,j}$ of the starting matrix $V_n$, we can use the assumptions (\ref{Ass0}, \ref{Ass1}, \ref{Ass2}) of Theorem \ref{th1}. Note that the fact that the variables $v_{i,j}$'s are i.i.d. allows us to group all the combinations which behave in the same way in the product of (\ref{arrang}). Let $\text{Part}(\cE)$ denote the set of all partitions of $\cE$, and set 

$$ \cE_{\gga}=\{\ovl{1}^1,\ldots,\ovl{k_1}^1\} \cup \cdots \cdots\cup \{\ovl{1}^p, \ldots,\ovl{k_p}^p\}.$$

For each partition $\ga$ in $\text{Part}(\cE)$, for each $x \in \cE$, we denote by $\ga(x)$   the  index of the class of $x$,  after having ordered the classes according to the order of their first element (for example, $\ga(\ovl{1}^1)=1$;  $\ga(\ovl{2}^1)=1$ if $\ovl{1}^1\stackrel{\ga}{\sim}\ovl{2}^1$ and  $\ga(\ovl{2}^1)=2$ if $\ovl{1}^1 \stackrel{\ga}{\nsim}\ovl{2}^1$).   
Therefore, we can write (\ref{arrang}) as two sums on the sets $\text{Part}(\cE)$, $\text{Part}(\cE_{\gga})$  introduced above. We get
 \[
 \E \left[ \prod_{l=1}^p \sum_{1\leq i\leq ns_l}(e^*_iX_n^{k_l} e_i-\E[e^*_iX_n^{k_l}e_i]) \right]  
 \]
 \begin{equation} \label{withpartitions}
= m^{-K} \sum_{\pi\in \text{Part}(\cE)}\; \sum_{\gga\in \text{Part}(\cE_{\gga})} A(n, \pi) B(m,\gga) \E \left[\prod_{l=1}^p \left(V_{l,\pi,\gga}-\E \, V_{l,\pi,\gga}\right)\right],
 \end{equation}
where: 
\begin{enumerate}
\item[$\bullet$] $V_{l,\pi,\gga}$ is defined this time with two functions $\pi$ and $\gga$ as shown in the previous paragraph for the general definition of $\ga$,
$$ V_{l,\pi,\gga} = v_{\pi(\ovl{0}^l),\gga(\ovl{1}^l)}v^*_{\gga(\ovl{1}^l),\pi(\ovl{1}^l)} \cdots  v_{\pi(\ovl{k_l -1}^l),\gga(\ovl{k_l}^l)}v^*_{\gga(\ovl{k_l}^l)\pi(\ovl{k_l}^l)}. 
$$
\item[$\bullet$] For each $\pi\in \text{Part}(\cE)$, $A(n,\pi)$ is the number of families of indices of  
 $$(i_{\ovl{0}^1}, \ldots,i_{\ovl{k_1}^1},i_{\ovl{0}^2}, \ldots,i_{\ovl{k_2}^2},\ldots\ldots,i_{\ovl{0}^p}, \ldots,i_{\ovl{k_p}^p})\in (1,\ldots,n)^{K}$$ whose level sets partition
is $\pi$ and that satisfies, for each $l= 1, \ldots,p$
\begin{align}   \label{cond0K}
1\leq \pi(\ovl{0}^l) = \pi(\ovl{k_l}^l) \leq ns_l.
\end{align}
\item[$\bullet$] For each $\gga \in \text{Part}(\cE_\gga)$, $B(m,\gga)$ is the number of families of indices of  
 $$(i_{\ovl{1}^1}, \ldots,i_{\ovl{k_1}^1},i_{\ovl{1}^2}, \ldots,i_{\ovl{k_2}^2},\ldots\ldots,i_{\ovl{1}^p}, \ldots,i_{\ovl{k_p}^p})\in (1,\ldots,m)^{K}$$ whose level sets partition
is $\gga$.
\end{enumerate}
$\,$

For any partitions $\pi \in \text{Part}(\cE)$ and $\gga \in \text{Part}(\cE_{\gga})$, let us define $G_{\pi,\gga}$ to be the graph with vertex set   
           $$V_{\pi,\gga} = \{\pi(x), \gga(x) \: ; \quad x \in \cE \:\, \mathrm{for} \:\, \pi, \quad x \in \cE_\gga \,\; \mathrm{for} \, \: \gga \},$$
 and edge set 
$$E_{\pi,\gga} =\llb \, \{\pi(\ovl{m-1}^l), \gga(\ovl{m}^l)\}, \{\gga(\ovl{m}^l),\pi(\ovl{m}^l)\}  \quad ; \; 1 \leq l \leq p, \;m\in\{1, \ldots,k_l\}\;\rrb.$$  
$\,$\\
\begin{center}
\tikz[scale=0.6]{
\draw[thick](-12,2)--(-1,2);\draw[thick,dashed](-1,2)--(1,2);
\draw[thick](-12,-2)--(-1,-2);\draw[thick,dashed](-1,-2)--(1,-2);
\draw[thick](1,2)--(13,2);
\draw[thick](1,-2)--(13,-2);
\draw (-12,2) node[left]{$\pi$};
\draw (-12,-2) node[left]{$\gga$};
\node (T1) at (-11,2) [above] {$0^1$};
\node (T2) at (-9,2) [above] {$1^1$};
\node (T3) at (-7,2) [above] {$2^1$};
\node (T4) at (-4,2) [above] {$(k_1-1)^1$};
\node (T5) at (-2,2) [above] {$k_1^1$};
\node (B1) at (-10,-2) [below] {$1^1$};
\node (B2) at (-8,-2) [below] {$2^1$};
\node (B3) at (-3,-2) [below] {$k_1^1$};
\node (TT1) at (2,2) [above] {$0^p$};
\node (TT2) at (4,2) [above] {$1^p$};
\node (TT3) at (6,2) [above] {$2^p$};
\node (TT4) at (8,2) [above] {$3^p$};
\node (TT5) at (10,2) [above] {$(k_p-1)^p$};
\node (TT6) at (12,2) [above] {$k_p^p$};
\node (BB1) at (3,-2) [below] {$1^p$};
\node (BB2) at (5,-2) [below] {$2^p$};
\node (BB3) at (7,-2) [below] {$3^p$};
\node (BB4) at (11,-2) [below] {$k_p^p$};
\draw[ultra thick](-11,2)--(-10,-2)(-9,2)--(-8,-2)(-4,2)--(-3,-2)(2,2)--(3,-2)(4,2)--(5,-2)(6,2)--(7,-2)(10,2)--(11,-2);
\draw[ultra thick,red](-9,2)--(-10,-2)(-7,2)--(-8,-2)(-2,2)--(-3,-2)(4,2)--(3,-2)(6,2)--(5,-2)(8,2)--(7,-2)(12,2)--(11,-2);
\node at (0,-4){\bf{Graph $G_{\pi,\gga}$}};
}
\end{center}
\vspace*{0.7cm}

For the term associated to a $(\pi, \gga) \in (\text{Part}(\cE) \times \text{Part}(\cE_\gga)) $ in  (\ref{withpartitions}) to be non zero, we need to have:
\begin{enumerate}
\item[(i)] for each $l=1, \ldots, p$, $\pi(\ovl{0}^l)=\pi(\ovl{k_l}^l)$, 
\item[(ii)] each edge of $G_{\pi,\gga}$ is visited at least twice by the union of the $p$
 paths \\
$(\pi(\ovl{0}^l),\gga(\ovl{1}^l),\pi(\ovl{1}^l), \ldots, \gga(\ovl{k_l}^l),\pi(\ovl{k_l}^l))$ , $\: l= 1, \ldots ,p\,$,
 \item[(iii)] for each $l=1, \ldots, p$, there exists $l' \neq l$ such that at least one edge of $G_{\pi,\gga}$ is visited by both paths  $(\pi(\ovl{0}^l),\gga(\ovl{1}^l),\pi(\ovl{1}^l), \ldots, \gga(\ovl{k_l}^l),\pi(\ovl{k_l}^l))$ and $(\pi(\ovl{0}^{l'}),\gga(\ovl{1}^{l'}),\pi(\ovl{1}^{l'}),\ldots, \gga(\ovl{k_{l'}}^{l'}),\pi(\ovl{k_{l'}}^{l'}))$. 
\end{enumerate}
Indeed, (i) is due to (\ref{cond0K}), (ii) is due to the fact that $v_{i,j}$'s are independent and centered and (iii) is due to the fact that the $v_{i,j}$'s are independent and that the variables $V_{l,\pi,\gga}-\E \, V_{l,\pi,\gga}$ are centered.

Now, let us define a function $s(\cdot)$ on the set $\cE$ in order to control the condition (\ref{cond0K}) in the following way: for each $l=1, \ldots,p$ and each $h=0, \ldots, k_l$, set  
$$ s({\ovl{h}^l})= \begin{cases} s_l & \textrm{if $h=0$ or $h=k_l$}\\
1&\textrm{otherwise,}
 \end{cases}$$
and 
\begin{equation} \label{spi}
s_\pi = \prod_{B \textrm{ bloc of } \pi} \min_{x \in B} s(x).
\end{equation}
 Then one can easily see that, as $n$ tends to infinity, 
$$
  A(n, \pi) \sim  s_\pi n^{|\pi|}\, ,
$$ 
 and 
$$
B(m, \gga) \sim  m^{|\gga|}\,,
$$
where $|\pi|$ denotes the number of vertices indexed by $\pi$ in the graph $G_{\pi,\gga}$, and $|\gga|$ denotes the number of vertices indexed by $\gga$ in $G_{\pi,\gga}$.

Therefore, for $(\pi,\gga)$ to have a non zero asymptotic contribution to (\ref{withpartitions}), we need the following  condition, in addition to (i), (ii) and (iii):  
\begin{enumerate}
\item[(iv)]  $  K \le |\pi|+|\gga|$.
\end{enumerate}
 
 Now, let us introduce this lemma which is the analogue of \cite[Lemma 2.1.34]{AGZ}. Its proof goes along the same lines as the proof of the former (see also \cite[Lemma 4.1]{BENA}).
 \begin{lemma}[Combinatorics]
 Let $(\pi, \gga) \in (\text{Part}(\cE) \times \text{Part}(\cE_\gga)) $ satisfy (i),(ii) and (iii). Then the number $C_{G}$ of connected components of $G_{\pi,\gga}$ is such that $\: C_G \le p/2 \: \, $ and 
 $$|V_{\pi,\gga}| = |\pi|+|\gga| \le C_G -\frac{p}{2}+ K. $$
 \end{lemma}
As a consequence, if $(\pi, \gga)$ also satisfies (iv), we have
\begin{enumerate}
\item[(a)] $C_G = p/2$,
\item[(b)] $p$ is even,
\item[(c)] $|V_{\pi,\gga}|= |\pi |+|\gga|= K$.
\end{enumerate}
Also note that by (ii), we have 
\begin{enumerate}
\item[(d)] $|E_{\pi,\gga}| \le K$,
\end{enumerate}
where $|E_{\pi,\gga}|$ denotes the number of edges of the graph $G_{\pi,\gga}$. Therefore, by (\ref{withpartitions}) and (c), we get
\[
\lim_{n \to \infty } \E \left[ \prod_{l=1}^p \sum_{1\leq i\leq ns_l}(e^*_iX_n^{k_l} e_i-\E[e^*_iX_n^{k_l}e_i]) \right]  = 
 \]
 \begin{equation} \label{with_lim}
 \sum_{\pi\in \text{Part}(\cE)}\; \sum_{\gga\in \text{Part}(\cE_{\gga})}\;  s_{\pi} \; \lim_{n \to \infty } \E \left[\prod_{l=1}^p \left(V_{l,\pi,\gga}-\E \, V_{l,\pi,\gga}\right)\right],
 \end{equation} 
where the sum is taken over the partitions $(\pi,\gga)$ which satisfy (i), (ii), (iii) and (iv) above, and such partitions also do  satisfy (a), (b), (c) and (d) above. \\

\noindent{\bf Case where $p$ is odd:} By (b), we know that when $p$ is odd, there is no couple of partitions $(\pi,\gga)$ satisfying the above conditions, hence
 \[
\lim_{n \to \infty } \E \left[ \prod_{l=1}^p \sum_{1\leq i\leq ns_l}(e^*_iX_n^{k_l} e_i-\E[e^*_iX_n^{k_l}e_i]) \right]  = 0.
 \]
\vspace*{0.5cm}
 
\noindent{\bf Case where $p=2$:} \label{p2} In this case, by (a) we know that for each couple of partitions $(\pi,\gga)$ satisfying (i), (ii), (iii) and (iv) above, the graph $G_{\pi,\gga}$ is connected. so that $|V_{\pi,\gga}|-1 \leq |E_{\pi,\gga}|$. Therefore, by (c) and (d) $|E_{\pi,\gga}| $ is either equal to $K$ or $K-1$:
\vspace*{0.4cm}
\begin{enumerate}
\item[$\bullet$] ${\bf |E_{\pi,\gga}|=K-1}$: In this case, the graph $G_{\pi,\gga}$ has exactly one more vertex than edges, hence it is a tree. As a consequence, the paths $(\pi(\ovl{0}^1),\gga(\ovl{1}^1),\pi(\ovl{1}^1),\gga(\ovl{k_1}^1),\pi(\ovl{k_1}^1))\,$  and $(\pi(\ovl{0}^2),\gga(\ovl{1}^2),\pi(\ovl{1}^2) ,\ldots,\gga(\ovl{k_2}^2),\pi(\ovl{k_2}^2)),$ which have the same beginning and ending vertices, satisfy the property that each visited edge is visited an even number of times. By an obvious cardinality argument, only one edge is visited more than twice, and it is visited four times (twice in each sense). The other edges are visited once in each sense. It follows that the expectation associated to a couple $(\pi,\gga)$ in (\ref{with_lim}) is equal to $\E[|v_{1,1}|^4]-1$.
\vspace*{0.3cm}
\item[$\bullet$]  ${\bf |E_{\pi,\gga}|=K}$: In this case, the graph $G_{\pi,\gga}$ has exactly the same number of vertices as edges, hence it is a bracelet. Therefore, by a cardinality argument again, the paths $(\pi(\ovl{0}^1),\gga(\ovl{1}^1),\pi(\ovl{1}^1)  ,\ldots, \gga(\ovl{k_1}^1),\pi(\ovl{k_1}^1))\,$  and $(\pi(\ovl{0}^2),\gga(\ovl{1}^2),\pi(\ovl{1}^2) ,\ldots, \gga(\ovl{k_2}^2),\pi(\ovl{k_2}^2))\,$ satisfy the property that they visit exactly twice of times each edge they visit (once in each sense). It follows that the expectation associated to a couple $(\pi,\gga)$ in (\ref{with_lim}) is equal to $1$.
\end{enumerate}

As a consequence, as $n$ tends to infinity, 
$$\E \left[\sum_{1\leq i\leq ns_1} (e^*_i X_n^{k_1} e_i- \E[e^*_i X_n^{k_1} e_i]) \times \sum_{1\leq i\leq ns_2} (e^*_i X_n^{k_2} e_i- \E[e^*_i X_n^{k_2} e_i]) \right]$$
 converges to a number that we shall denote by  
 \begin{align} \label{cov}
\Cov_{s_1,s_2}(k_1,k_2),
 \end{align}
 which  depends on the distribution of the $v_{i,j}$'s only through $\lim_{n \to \infty} \E[|v_{1,1}|^4]$. 
\vspace*{0.3cm}

\noindent{\bf Case where $p$ is $>2$ and even:} By (a) above, for each couple of partitions $(\pi,\gga)$ satisfying (i), (ii), (iii) and (iv),  $G_{\pi,\gga}$ has exactly $p/2$ connected components. By (iii), each one of them contains the support of exactly two of the $p$ paths
$$(\pi(\ovl{0}^l),\gga(\ovl{1}^l),\pi(\ovl{1}^l) \ldots \ldots, \gga(\ovl{k_l}^l),\pi(\ovl{k_l}^l))\, \quad \quad \quad(l= 1, \ldots ,p).$$
To join every two paths having the same support in the expectation of (\ref{with_lim}), let us define $\gs_{\pi,\gga}$ to be the  matching (i.e. a permutation all of whose cycles have length two) of $\{1, \ldots, p\}$ such that for all $ l=1 , \ldots, p \,$, the paths with indices $l$ and $\gs_{\pi,\gga}(l)$ are supported by the same connected component of  $G_{\pi,\gga}$. 

We shall now partition the sum of (\ref{with_lim}) according to the value of the matching $\gs_{\pi,\gga}$ defined by $(\pi,\gga)$. We get 
\begin{equation}
\lim_{n \to \infty } \E \left[ \prod_{l=1}^p \sum_{1\leq i\leq ns_l}(e^*_iX_n^{k_l} e_i-\E[e^*_iX_n^{k_l}e_i]) \right]  = 
\end{equation}
 \[ \sum_{\gs} \sum_{\pi,\gga}  s_{\pi} \; \lim_{n \to \infty } \E \left[\prod_{l=1}^p \left(V_{l,\pi,\gga}-\E \, V_{l,\pi,\gga}\right)\right] \mathbf{1}_{\{\gs_{\pi,\gga}=\gs\}},
  \]
where the first sum is over the matchings $\gs$ of $\{1, \ldots, p\}$  and the second sum is over the couples of partitions $(\pi,\gga)$ satisfying (i), (ii), (iii) and (iv).\\

Note that for each matching $\gs$ of $\{1, \ldots, p\}$, the set of  couples of  partitions $(\pi,\gga) \in (\text{Part}(\cE) \times \text{Part}(\cE_\gga)) $ such that $\gs_{\pi,\gga}=\gs$ can be identified with the Cartesian product, indexed by the set  of cycles $\{l,l'\}$ of $\gs$, of the set of couples of partitions $(\ovl{\pi},\ovl{\gga}) \in (\text{Part}(\cE^{l,l'}) \times \text{Part}(\cE_\gga^{l,l'}))$ such that
$$ \cE^{l,l'}  = \{\ovl{0}^l, \ldots,\ovl{k_l}^l\} \, \cup \, \{\ovl{0}^{l'}, \ldots,\ovl{k_l}^{l'}\} \quad \quad \quad (\textrm{subset of} \; \cE)
$$
and 
$$ \cE_{\gga}^{l,l'}=\{\ovl{1}^l, \ldots,\ovl{k_l}^l\} \, \cup \, \{\ovl{1}^{l'}, \ldots,\ovl{k_l}^{l'}\} \quad \quad \quad (\textrm{subset of} \; \cE_{\gga})
$$
satisfying the following conditions
\begin{enumerate}
\item[(i')]$\ovl{\pi}(\ovl{0}^l)=\ovl{\pi}(\ovl{k_l}^l)$  and $\ovl{\pi}(\ovl{0}^{l'})=\ovl{\pi}(\ovl{k_{l'}}^{l'})$,
\item[(ii')] each edge of the graph $G_{\ovl{\pi},\ovl{\gga}}$  is visited at least twice by the union of its two paths indexed by the corresponding $l$ and $l'$.
\item[(iii')]at least one edge of $G_{\ovl{\pi},\ovl{\gga}}$ is visited by both previous paths,
\item[(iv')] $k_l+k_{l'} \leq  |\ovl{\pi}|+|\ovl{\gga}|$.
\end{enumerate} 
Moreover, one can see that the factor $s_{\pi}$ factorizes along the connected components of $G_{\pi,\gga}$, and by the independence of the random variables $v_{i,j}$'s, the expectation 
$$ \E \left[ \prod_{l=1}^p \sum_{1\leq i\leq ns_l}(e^*_iX_n^{k_l} e_i-\E[e^*_iX_n^{k_l}e_i]) \right]$$
also factorizes along the connected components of $G_{\pi,\gga}$. Therefore, we get
\begin{equation}
\lim_{n \to \infty } \E \left[ \prod_{l=1}^p \sum_{1\leq i\leq ns_l}(e^*_iX_n^{k_l} e_i-\E[e^*_iX_n^{k_l}e_i]) \right]  =  
\end{equation}
 \[ \sum_{\gs} \prod_{l,l'}  \; \lim_{n \to \infty } \E \left[\sum_{1\leq i\leq ns_l}(e^*_iX_n^{k_l} e_i-\E[e^*_iX_n^{k_l}e_i]) \times \sum_{1\leq i\leq ns_{l'}}(e^*_iX_n^{k_{l'}} e_i-\E[e^*_iX_n^{k_{l'}}e_i])\right]
  \]
where the sum is over the matchings $\gs$ of $\{1, \ldots, p\}$ and for each such $\gs$, the product is over the cycles $\{l,l'\}$ of $\gs$.

By the previous definition of $\Cov_{s_1,s_2}(k_1,k_2)$ in $(\ref{cov})$, we get
\[
\lim_{n \to \infty } \E [ \prod_{l=1}^p \sum_{1\leq i\leq ns_l}(e^*_iX_n^{k_l} e_i-\E[e^*_iX_n^{k_l}e_i]) ]  =  \sum_{\gs \textrm{ matching}} \; \prod_{\{l,l'\} \textrm{cycle of } \gs} \Cov_{s_l,s_{l'}}(k_l,k_{l'}).
\] 
By Wick's formula and Equation (\ref{centrage}), we have proved the first part of Proposition (\ref{mainprop}). \\

We finish the proof by this last step.
\vspace{0.3cm}

\noindent{\bf Computation of $\Cov_{s_1,s_2}(1,1)$:} In this case, we have $p=2,$ and $K=k_1+k_2=2$. Therefore, by (a), (c) and (d),
\begin{enumerate}
\item[$\bullet$] $G_{\pi,\gga}$ is connected, 
\item[$\bullet$] $V{\pi,\gga}=K=2$,
\item[$\bullet$] $E_{\pi,\gga} \leq 2$.
\end{enumerate} 
With two vertices, there is exactly one tree and zero bracelet. Thus, we have by the paragraph devoted to the case $p=2$,
$$ \Cov_{s_1,s_2}(1,1) = \sum_{(\pi, \gga), G_{\pi, \gga} \textrm{ is a tree}} s_{\pi} \: (\lim_{n \to \infty} \E[|v_{1,1}|^4]-1)$$ 
 For this tree, there are two associated couples of partitions $(\pi,\gga)$: 

\begin{center}
\tikz[scale=0.6]{
\draw[very thick](0,-3.5)--(0,2.5);
\draw[thin](-10,2)--(-0.5,2);\draw (0.5,2)--(10,2);
\draw[thin](-10,-2)--(-0.5,-2);\draw (0.5,-2)--(10,-2);
\draw (-10,2) node[left]{$\pi$};
\draw (-10,-2) node[left]{$\gga$};
\node (T1) at ( -9,2) [NN] {};
\node (T2) at ( -7,2) [NN] {};
\node (T3) at ( -4,2) [NN] {};
\node (T4) at ( -2,2) [NN] {};
\node (B1) at ( -8,-2) [NR] {};
\node (B2) at ( -3,-2) [NR] {};
\node (T5) at ( 2,2) [NN] {};
\node (T6) at ( 4,2) [NN] {};
\node (T7) at ( 7,2) [NR] {};
\node (T8) at ( 9,2) [NR] {};
\node (B3) at ( 3,-2) [NR] {};
\node (B4) at ( 8,-2) [NN] {};
\draw[ultra thick](T1)--(B1);
\draw[ultra thick,red](B1)--(T2);
\draw[ultra thick](T3)--(B2);
\draw[ultra thick,red](B2)--(T4);
\draw[ultra thick](T5)--(B3);
\draw[ultra thick,red](B3)--(T6);
\draw[ultra thick](T7)--(B4);
\draw[ultra thick,red](B4)--(T8);
\node  at (-5.5,-3){\bf{Case 1}};
\node at (6,-3){\bf{Case 2}};
}
\end{center}
$\,$

 {\bf Case (1)}: The partition $\pi$ is defined by
$$ \pi(\cE)= \{\{\ovl{0}^1,\ovl{1}^1,\ovl{0}^2,\ovl{1}^2\}\},$$
 hence by (\ref{spi}), $s_{\pi} = \min\{s_1,s_2\}$.\\

 {\bf Case (2)}: In this case, the partition $\pi$ is defined by
$$\pi(\cE)= \{\{\ovl{0}^1,\ovl{1}^1\} , \{\ovl{0}^2,\ovl{1}^2\}\},$$
 hence, $s_{\pi} = s_1 s_2$.\\

  As a consequence, $$\Cov_{s_1,s_2}(1,1)= (\lim_{n \to \infty} \E[|v_{1,1}|^4]-1) (\min\{s_1,s_2\}+s_1 s_2).$$
Now, by Equation (\ref{centrage}), we get
\[
\begin{split}
& \lim_{n \to \infty}\E \left[ \sum_{1\leq i\leq ns_1}(e^*_iX_n e_i-\frac{1}{n} \text{Tr}(X_n))\times \sum_{1\leq i\leq ns_2}(e^*_iX_n e_i-\frac{1}{n} \text{Tr}(X_n))\right]   \\
&= \Cov_{s_1,s_2}(1,1)- s_2 \, \Cov_{s_1,1}(1,1)-s_1 \, \Cov_{1,s_2}(1,1)+s_1s_2 \, \Cov_{1,1}(1,1)  \\
&= (\lim_{n \to \infty} \E[|v_{1,1}|^4]-1)\left[ (\min\{s_1,s_2\}+ s_1s_2)-( 2 s_1 s_2  )-( 2 s_1 s_2  )+(2 s_1 s_2)\right] \\
&= (\lim_{n \to \infty} \E[|v_{1,1}|^4]-1)(\min\{s_1,s_2\}- s_1s_2), 
 \end{split}
 \]
which concludes the proof of Proposition \ref{mainprop}. {\hfill $\square$}
\section{Tightness and Convergence in the Skorokhod topology} \label{tightness}

For Wigner matrices, Benaych-Georges proved that the bivariate process $B^n$ converges in distribution, for the Skorokhod topology in $D[0,1]^2$, to the  bivariate Brownian bridge under several assumptions on the atom distribution: absolute continuity, moments of all orders and matching with a GUE/GOE matrix up to order 10 on the diagonal and up to order 12 off the diagonal. In order to prove this convergence, he used some ideas developed by Tao and Vu  in \cite{Tao_Vu}, especially, the "Four Moment Theorem for eigenvectors of  Wigner matrices", see \cite[theorem 8]{Tao_Vu}. \\    % that we heavily believe

To our knowledge, such a theorem is not yet available for the case of sample covariance matrices. We formulate the statement in Hypothesis \ref{hypo} below. If this is indeed the case, convergence of the process $B^n$ for the Skorokhod topology in $D[0,1]^2$ will be also verified in our case. For proving this, we will follow closely the strategy of Benaych-Georges.

\begin{definition}[Matching moments] Let $k \geq 1$. Two random matrices $V_n=(v_{i,j})_{i \leq n,\, j\leq m}$, $V^{'}_n=(v^{'}_{i,j})_{i \leq n,\, j\leq m} $ are said to match up to order $k$, if one has 
  $$ \E Re(v_{i,j})^a Im(v_{i,j})^b =   \E Re(v^{'}_{i,j})^a Im(v^{'}_{i,j})^b $$
whenever $a,b \geq 0$, $1\leq i \leq n$ and $ 1 \leq j \leq m $ are integers such that $\: a+b \leq k$.
\end{definition}
Before stating the theorem of convergence of our process $B^n$ for the Skorokhod topology in $D[0,1]^2$, let us give the following  hypothesis that we will need: 

 \begin{hyp}[Matching theorem for eigenvectors] \label{hypo} 
We suppose that if the matrix $V_n$ matches a ($n \times m$) - Gaussian matrix $M_n$ (i.e. matrix whose elements are independent standard Gaussian variables) to order $l \geq 4$, then, for any fixed positive integer $k$ and polynomial function $G$ on $\C^k$, there exists a certain constant C independent of $n$ such that  
\begin{align}
|\E[G(n\,u_{p_1,i_1}\ovl{u}_{q_1,i_1}, \ldots,n\,u_{p_k,i_k}\ovl{u}_{q_k,i_k})]-\E[G(n\,u'_{p_1,i_1}\ovl{u}'_{q_1,i_1}, \ldots,n\,u'_{p_k,i_k}\ovl{u}'_{q_k,i_k})]\;|\le Cn^{2-\frac{l}{2}}.
\end{align}
whenever  $(i_1,p_1,q_1),\ldots,(i_k,p_k,q_k)$ is a collection of indices in $\{1,\ldots,n\}^3$, $U_n=(u_{i,j})_{1 \leq i,j \leq n}$ is the eigenmatrix of the sample covariance matrix $X_n=(1/m) V_n V^{*}_n $ and  $U'_n=(u'_{i,j})_{1 \leq i,j \leq n}$ is the eigenmatrix of the Laguerre matrix $L_n=\frac{1}{m}M_n M^{*}_n$.
 \end{hyp}
\begin{theorem}[Convergence in the Skorokhod topology] \label{thm2} For the sample covariance matrix $X_n=(1/m) V_n V^{*}_n$ defined as in Section \ref{sec2}, suppose that
\begin{enumerate}
\item[(i)] The distribution of the entries of $\, (V_n)$ are absolutely continuous with respect to the Lebesgue measure.
\item[(ii)]  $ \forall \, k \geq 0, \: \; \sup_n \, \E|v_{1,1}^{(n)}|^{k}< \infty.$
\item[(iii)] $(V_n)$ matches a ($n \times m$) - Gaussian matrix $M_n$ up to order $l$, and that Hypothesis \ref{hypo} is satisfied.
\end{enumerate}  
Then, for $l=12$, the bivariate process $B^n$ converges in distribution, for the Skorokhod topology in $D([0,1]^2)$,  to the bivariate Brownian bridge. 
\end{theorem}
\begin{rem}[Comments on the assumptions of Theorem \ref{thm2}] These assumptions might not to be optimal, especially the continuity one and matching up to order 12. We hope to prove this theorem under Assumption (iii) for l=4 instead of l=12. 
\end{rem}
 Proving Theorem \ref{thm2} consists to prove the following lemma of tightness and uniqueness of the accumulation point argument.
 \begin{lemma}[Tightness argument]  \label{lem_tight}
 Under Assumptions of Theorem \ref{thm2}, the sequence (distribution $(B^n))_{n\geq 1}$ is $C$-tight, i.e.  is tight and has only one  $C([0,1]^2)$-supported accumulation point.
 \end{lemma}
\subsubsection{Proof of Theorem \ref{thm2}}
 Note that Theorem \ref{thm2} allows us to show that for all $0\leq s<s'\leq 1$ and $0\leq t<t'\leq 1$, the sequence of random variables 
 $$\frac{1}{\sqrt{(s'-s)(t'-t)}}\sum_{\substack{ns < i \leq ns' \\ nt<j\leq nt'}}(|u_{i,j}|^2-1/n)$$ 
admits a limit in distribution as $n\to\infty$, hence is bounded in probability (in the sense of \cite[Def. 1.1: $  \lim_{C \to \infty} \liminf_{n \to \infty} \bbP(|X_n| \leq C)=1$]{Tao_Vu_2}). In the next proposition, we improve these assertions by making them uniform on $s,s',t,t',i,j$ and upgrading them to the $L^2$ and $L^4$ levels. This proposition is almost sufficient to apply the tightness argument of the (distribution $(B^n))_{n\geq 1}$.   
\begin{prop}[Control of Jumps]\label{propL2L4} Suppose that Assumptions (i), (ii) and (iii) for $l=4$  (resp. $l=8$) are satisfied. 
 Then as $n \to \infty$, the sequence   
 \begin{align} \label{eqp}
 n|u_{i,j}|^2-1 \qquad \textrm{( resp.  } \qquad\frac{1}{\sqrt{(s'-s)(t'-t)}}\sum_{\substack{ns<i\leq ns'\\ nt<j\leq nt'}}(|u_{i,j}|^2-1/n)\textrm{)}
 \end{align} 
is bounded for the  $L^4$ (resp. $L^2$) norm, uniformly in $i,j$ (resp. $s<s',t<t'$).
 \end{prop}
The proof of Proposition \ref{propL2L4} goes along the same lines as the proof given by Benaych-Georges in \cite[section 4.4]{BENA}. Indeed, Hypothesis \ref{hypo} "matching with Gaussian matrix" allows us to work with the entries of a Haar-distributed matrix instead of the entries of the eigenmatrix of the sample covariance matrix $X_n$. Note only that, if the second term of (\ref{eqp}) have been bounded for $L^{2+\epsilon} $ instead of $L^2$, Assumption (iii) for $l=8$ would have been enough to prove the convergence of $B^n$ in distribution, for the Skorokhod topology in $D([0,1]^2)$,  to the bivariate Brownian bridge.    \\

Now, to prove Lemma  \ref{lem_tight}, we give the following proposition, which is the obvious multidimensional generalization of Proposition 3.26 of \cite[Chapter VI]{Jac_Shi}: \\

For $f\in D([0,1]^2)$ and  $(s_0,t_0)\in [0,1]^2$, we define $\Delta_{s_0,t_0} f$ to be the  "maximal jump" of $f$ at $(s_0,t_0)$, i.e. 
$$\Delta_{s_0,t_0} f:=\max_{\diamond, \diamond'\in\{<,\geq\}}\left|f(s_0,t_0)-
 \lim_{\substack{s \to_{\diamond} s_0\\ t \to_{\diamond'} t_0}}f(s,t)\right|.$$
 \begin{prop}[C-Tightness]
If the sequence $($distribution $(B^n))_{n\geq 1}$ is tight and satisfies 
\begin{align}\label{delta}
\forall \epsilon >0,\qquad \bbP(\sup _{(s,t)\in [0,1]^2}\Delta_{s,t} B^n >\epsilon)\rightarrow_{n \to \infty}  0,
\end{align}
then the sequence $($distribution $(B^n))_{n\geq 1}$ is $C$-tight, i.e.  is tight and has only one $C([0,1]^2)$-supported accumulation point.
 \end{prop}
So to prove Lemma \ref{lem_tight}, let us first prove that the sequence $($ distribution $(B^n))_{n\ge 1}$ is tight.  For this, we follow closely the proof of Benaych-Georges.
 
Note that   the process $B^n$ vanishes at the border of $[0,1]^2$.  So according to \cite[Th. 3]{Bic_Wic} and to Cauchy-Schwartz inequality,  
it suffices to prove that there exists $C<\infty$  such that for $n$ large enough,  for all   $s<s', \, t<t'\in [0,1]$, 
$$\E[\{\sum_{ns<i\leq ns'}\sum_{nt< j\leq nt'}(|u_{i,j}|^2-1/{n})\}^{4}]\leq C(s'-s)^{2}(t'-t)^{2}.$$
As in the proof of Proposition \ref{propL2L4}, one can suppose that the $u_{i,j}$'s are the entries of a Haar-distributed matrix. But in this case, the job has already   been done in \cite{Don_Rou}: the unitary case is treated in  Section 3.4.1 (see specifically  Equation (3.25)) and the orthogonal case is treated, more elliptically,  in Section 4.5.  
%(to   recover the details of the proof, join Equations (3.26), (4.5) and (4.17)).
 
Let us now prove (\ref{delta}). Note that $\sup _{(s,t)\in [0,1]^2}\Delta_{s,t} B^n =\max_{1\leq i,j\leq n} ||u_{i,j}|^2-1/n|$. As a consequence, by the union bound, it suffices to prove that for each $\epsilon >0$, there exists $C <\infty$ independent of $i,j$ and $n$ such that for all $i,j$, 
$$ \bbP(||u_{i,j}|^2-1/n |>\epsilon)\leq Cn^{-4},$$
 which follows from Chebyshev's inequality and Proposition \ref{propL2L4}. \hfill$\square$

\bibliography{bibliographie}
\bibliographystyle{plain}

\end{document}